\documentclass[a4paper,12pt]{amsart}

\usepackage[T1]{fontenc}
\usepackage{amsmath,amssymb,amscd,amsfonts,amsthm,mathrsfs,eucal,dsfont,verbatim}
\RequirePackage{doi}
\usepackage{lscape,enumerate}
\usepackage{color,soul}

\usepackage[margin=2in]{geometry}
\usepackage[square,numbers,sort&compress]{natbib}
\usepackage{marginnote}

\numberwithin{equation}{section}

\newcommand{\real}{\mathds{R}}

\newcommand{\Ee}{\mathds E}
\newcommand{\Pp}{\mathds P}
\newcommand{\I}{\mathds 1}

\newcommand{\Aa}{\mathcal{A}}

\newcommand{\LL}{\mathcal{L}}
\newcommand{\Bb}{\mathcal{B}}
\def\1{1\!\!\hbox{{\rm I}}}

\newcommand{\N}{\mathcal{N}}

\newcommand{\rd}{{\mathds{R}^d}}

\newcommand{\sign}{\operatorname{sign}}

\usepackage{color}

\evensidemargin
-1cm

\usepackage[T2A]{fontenc}
\usepackage[english]{babel}

\numberwithin{equation}{section}
\theoremstyle{plain}
\newtheorem{theorem}{Theorem}
\newtheorem{lemma}{Lemma}

\theoremstyle{definition}

\newtheorem{example}{Example}

\numberwithin{theorem}{section}
\numberwithin{definition}{section}
\numberwithin{corollary}{section}
\numberwithin{remark}{section}
\numberwithin{lemma}{section}

\counterwithout{equation}{section}
\counterwithout{theorem}{section}
\counterwithout{lemma}{section}

\begin{document}

	\title{On ergodic property of   the solution to a L\'evy-driven SDE}

	\date{}

	\author[V. Knopova]{Victoria Knopova}
	\address[V. Knopova]{Kiev T. Shevchenko University\\ Department of Mechanics and Mathematics\\ Acad. Glushkov Ave., 02000,  Kiev, Ukraine\\ \texttt{Email: $vicknopova$@$knu.ua$ (corresponding author)}}

	\author[Y. Mokanu]{Yana Mokanu}
	\address[Y. Mokanu]{I. Sikorsky Kyiv Polytechnic Institute\\ Dep. of Mathematical Analysis and Probability Theory\\ Beresteyskyi Ave. 37, 03056 Kyiv, Ukraine\\
	\texttt{Email:  $yana.mokanu3075$@$gmail.com$} }

	\begin{abstract}
		In this paper, we investigate  ergodicity in total variation of the process $X_t$, related to  a  L\'evy-driven stochastic differential equation with unbounded coefficients, and describe the speed of convergence to the respective invariant measure. Some examples are provided.
	\end{abstract}



	\maketitle 
	
	\section{Problem Statement}\label{PS}
	
	Let $Z$ be a L\'evy process on $\rd$. Denote by $\N(du,ds)$ its Poisson random measure, by $\nu(du) ds$ its compensator, and consider a  L\'evy-driven SDE (in the integral form): 
	\begin{equation}\label{SDE}
		\begin{split}
			X(t) &= X(0) + \int_0^t \ell(X(s))ds
			+ \int_0^t \int_{u\neq 0} k(X(s-),u)\I_{|k(X(s-),u)|> 1}\N(du,ds)\\
			& \quad + \int_0^t \int_{u\neq 0} k(X(s-),u)\I_{|k(X(s-),u)|\leq 1}\left(\N(du,ds)- \nu(du)ds)\right), 
		\end{split}
	\end{equation}
	where $\ell(\cdot): \rd \to \rd$,  $k(\cdot,\cdot): \rd\times \rd \to \rd$ are measurable functions. 
	
	We assume that there exists  (weak or strong)  solution to \eqref{SDE}, which is  a strong Markov process; the respective conditions could be found in \cite{GS82} or \cite{IW89}. 
	
	It can be shown that for any function from the space  $C^2_0(\rd)$  of twice continuously differentiable functions with compact support,  the process 
	\begin{equation}
		M_t := f(X_t)- \int_0^t L f(X_s)ds, 
	\end{equation}
	where 
	\begin{equation}\label{L01}
		\begin{split}
			Lf(x) &=  \ell(x)\cdot \nabla f(x)\\
			& \quad + \int_{\rd \backslash \{0\}} \left(f(x+k(x,u)) - f(x) -\nabla f(x) \cdot k(x,u) \I_{|k(x,u)|\leq 1} \right) \nu(du)
		\end{split}
	\end{equation}
	is a martingale (cf. \cite{GS82}). The domain of $L$ can be extended (see e.g. the proof of the Ito formula (cf. \cite{IW89}, \cite{Si05})  to the space $C^2(\rd)$. 
	We call $(\LL, C^2(\rd))$ the extended generator of the semigroup $(P_t)_{t\geq 0}$, $P_tf(x) = \Ee^x f(X_t)$.

	In this paper we investigate  ergodicity in total variation of the process $X_t$ related to \eqref{L01} and  describe the speed of convergence to the  invariant measure. The usual path for such investigation is the following: 
	if  the transition probability kernel satisfies certain irreducibility assumption (e.g. the local Dobrushin condition) and if the  extended generator  $(\LL, C^2(\rd))$  satisfies the Lyapunov-type   condition  with some Lyapunov function,  then one can reduce the continuous time  setting to the discrete one (cf., e.g.,   \cite{Ku18}) and apply ergodic theorems  to  the skeleton chain  $X^h := \{ X_{nh}, \, n\in \mathbb{Z}_+\}$.  Under  irreducibility assumption,  ergodicity of the skeleton chain implies the ergodicity of the initial Markov process with the same convergence rates.

	The existence of an invariant measure and the speed of convergence to it are important, particularly in simulating stochastic processes. A prime example is their role in investigating the existence of a stationary state in a stochastic system and in simulating functionals related to the stationary distribution. Information about the speed of convergence to the stationary distribution allows us to evaluate how accurately the stationary solution is approximated by simulating the original  stochastic differential equation.

	This paper is in a sense a continuation of the investigation started in  \cite{KM24}, where the one-dimensional case was studied in the context of L\'evy-type processes  and the L\'evy-type kernel was assumed to be bounded in $x$.  In the  multi-dimensional case the situation is much more complicated due to the geometry of the support of the  L\'evy measure $\nu(\cdot)$.
	
	Loosely speaking, ergodicity of the solution to \eqref{SDE} can be induced by the drift as well as by the noise.
	It is a kind of a ``general knowledge'' that in a  pure diffusion setting  one can get the Lyapunov condition for $L$ due to the fast growth of the off-diagonal coefficients; in this  case we have the  noise-induced ergodicity.  Such examples can be deduced from \cite{Kh12}, which is a classical reference on ergodic Markov processes.  Motivated by such examples, we try to answer the question, which assumptions on the support of the L\'evy-type  kernel provide the Lyapunov-type condition  for the  L\'evy-type operators given by \eqref{L01}. We are interested in the  ``noise-induced'' ergodicity  as well as in its interplay with the ``drift-induced'' one. 
	With some efforts, this result can  be  reformulated for the situation,  when we have the varying intensity of the process rather than varying jump size.
	
	Let us give a brief overview of the existing results. While for diffusion processes the literature is rather rich, for jump-diffusion and general  L\'evy-type processes the literature is not so extensive.  The pioneering works devoted to the topic belong to Skorokhod; in \cite{Sk87} a long-time behaviour of general Markov processes and solutions to SDEs was discussed, as well as stabilization properties of solutions to systems of SDEs.
	Ergodicity of a solution to the jump-diffusion SDE  was studied in Wee \cite{We99}, where the conditions are given in the integral form. Masuda  in  \cite{Ma04} investigated  a L\'evy-driven Ornstein-Uhlenbeck process, in particular, the existence of a smooth density, ergodicity and convergence rates; later  in \cite{Ma07} a general case of jump-diffusion with compound Poisson jump measure was investigated.  Exponential ergodicity of  a solution  to a L\'evy-driven  SDE  was  investigated in Kulik  \cite{Ku09}, as well as the assumptions, which imply the Doeblin condition.  See also Wang~\cite{Wa08} for  conditions for a 1-dimensional L\'evy-driven SDE.  In Fort, Roberts \cite{FR05}  and Douc, Fort, Guillin \cite{DFG09}  sufficient conditions for ergodicity are given in the form of so-called $(f,r)$-modulated moments. 
	For a general L\'evy-type process   conditions ensuring  ergodicity  were  studied by Sandri\'c in the series of papers \cite{Sa13}, \cite{Sa14}, \cite{Sa16}, in particular, in  \cite{Sa16}, where the conditions are provided in a (very general)  integral  form.

	The paper is organized as follows.  In Section~\ref{Prel},   we recall the notion of ergodicity and  discuss the assumptions.  In Section~\ref{Set},  we formulate our results, in particular,   we provide the sufficient conditions  under which the  process is  ergodic and  derive  the speed of convergence to the invariant  measure. Proofs are given in Section~\ref{Proofs}. We illustrate our results by examples in Section~\ref{Examples}.

	\section{Preliminaries }\label{Prel}

	Recall that a Markov process $(X_t)_{t\geq 0}$  on $\rd$  with a transition probability kernel $P_t(x,dy)$ is called (strongly)
	\emph{ergodic} if  there exists an  invariant probability measure $\pi(\cdot)$ such that 
	\begin{equation}\label{inv1}
		\lim_{t\to \infty} \| P_t(x,\cdot)- \pi(\cdot)\|_{TV} =0
	\end{equation}
	for any $x\in \real^d$.
	Here $\|\cdot\|_{TV}$ is the total variation norm. Under additional assumptions one can find the  speed of convergence to the invariant measure. We say, that the process   is  (non-uniformly) $\psi$-ergodic, if there exists a function $\psi\geq 0$, $\psi(t) \to 0$ as $t\to\infty$, such that 
	\begin{equation}\label{psierg}
		\|P_t(x,\cdot )- \pi(\cdot)\|_{TV}\leq C(x) \psi(t), \quad t\geq 1, 
	\end{equation}
	with some $C \colon \rd \to [0,\infty)$. 
	
	As we  already mentioned in Section~\ref{PS}, one of the possible ways to show the ergodic property relies on the following milestones: reduction of the initial Markov process to the skeleton chain,  (local)  Dobrushin condition   and the  Lyapunov-type  condition  for  this chain, and finally the ``return route''  to the initial process. 
	Let us recall these conditions. 
	
	We say that a  Markov chain $X$  with a transition kernel $\Pp(x,dy)$  satisfies the local Dobrushin condition (cf. \cite{Ku18}) on a measurable set $B\subset \rd\times \rd$  if 
	\begin{equation}\label{Dor1}
		\sup_{x,y\in B} \|\Pp(x,\cdot)- \Pp(y,\cdot)\|_{TV} <2. 
	\end{equation}
	The sufficient condition for \eqref{Dor1} is provided by the existence of an absolutely continuous component of the transition kernel  $\Pp(x,dy)$ of the chain, such that $x\mapsto \Pp(x,\cdot)$ is continuous in $L_1$,   see \cite{Ku18}. This condition is clearly satisfied if the initial Markov process has a strictly positive transition probability density $p_t(x,\cdot)$,  such that $x\mapsto p_t(x,\cdot)$ is continuous in $x$ in $L_1$. 
	On the language of semigroups this condition provides that the respective semigroup $(P_t)_{t\geq 0}$, $P_t f(x):= \Ee^x f(X_t)$, possesses the $C_b$-Feller  property.  Dobrushin condition ensures, that the processes with different starting points are attracted at infinity by the same measure; in other words, this condition ensures that two processes can be coupled at infinity. 
	
	The second   condition we need  is the \emph{Lyapunov-type} inequality.   This inequality  allows to establish moment estimates, responsible for (non-uniform) ergodic rates of convergence. Let us recall this condition  for the skeleton chain.
	
	Suppose that there exist 
	
	\begin{itemize}
		\item   a  norm-like function $V: \, \rd\to [1,\infty)$ (i.e. $V(x)\to \infty$ as $|x|\to  \infty$),   bounded on a compact $K$; 
		
		\item  a function $f: \, [1,\infty) \to (0,\infty)$, that  admits a non-negative increasing and concave extension to $[0,\infty)$;
		
		\item  a constant $C>0$,
	\end{itemize}
	such that the following  relations  hold: 
	\begin{equation}\label{Lyap1}
		\Ee_x V(X_1) - V(x) \leq - f (V(x))+ C, \quad x\in \rd, 
	\end{equation}
	\begin{equation}\label{Lyap11}
		f\left(1+ \inf_{x\notin K} V(x)\right) > 2C.  
	\end{equation}
	Inequality \eqref{Lyap1} is called the \emph{Lyapunov-type  inequality}, and the respective  function $V$ is called the \emph{Lyapunov function}.  
	
	Since we start not with  a Markov chain but with  a Markov process, it is more convenient to check a condition, which ensures that the respective skeleton process possesses the Lyapunov-type condition. The required  condition is formulated in terms of the extended  generator $\LL$ (because Lyapunov functions are usually unbounded, we need to extend the domain of $L$ to a bigger space):
	\begin{equation}\label{Lyap2}
		\LL  V(x) \leq - f(V(x))  + C, 
	\end{equation}
	where $V$ and $f$ have the same meaning as above.   
	
	To summarize,  if we manage to prove \eqref{Lyap2}, then by \cite{Ku18} inequality \eqref{Lyap1} holds true for the skeleton chain. Lyapunov  and Dobrushin conditions  allow to prove  (cf. \cite{Ku18}) that  the chain $X$ is ergodic  and to find the (non-uniform) ergodic rate.   This theorem together with the Dobrushin condition ensures the ergodicity of the chain and yields the rate of convergence to the ergodic ditribution. Moreover, the same convergence rate will still hold for the initial Markov process, see the discussion in  \cite{Ku18}, also \cite{Sa16}. The speed of  convergence to $\pi$  in \eqref{psierg} as $t\to \infty$ is (cf. \cite{Ku18})
	\begin{equation}\label{psi}
		\psi (t) := \left(\frac{1}{f(F^{-1}(\gamma t))}\right)^\delta, \qquad F(t) := \int\limits_1^t \frac{dw}{f(w)}, \quad t \in (1, \infty),
	\end{equation}
	where $\gamma,\delta\in (0,1)$. Moreover, 
	\begin{equation*}
		\int_{\rd} f(V(x)) \pi (dx) < \infty. 
	\end{equation*}
	
	The approach described above is not the only possible one; for the discussion of other possible ways of showing ergodicity we refer to \cite{Ku18}.

	\section{Settings and results}\label{Set}

	Our standing assumptions in this paper  would be  
	
	\begin{align}\label{A1}
		\tag{\bfseries{A1}}
		&\text{There exists   (weak or strong)   a  solution $\{X_t, \,t>0\}$ to \eqref{SDE}}\\ \label{A2}\tag{\bfseries{A2}}
		&
		\text{The  skeleton chain  $X^h := \{ X_{nh}, \, n\in \mathbb{Z}_+\}$ satisfies  the Dobrushin condition.}	
	\end{align}
	Condition  \eqref{A2} provides the required irreducibility of  the   skeleton process.

	As usual, we denote by   $|u|$ the norm of a vector $u\in \rd$. 
	For  $\phi: \rd\to \real_+$, $\phi\in C^2(\rd)$, such that $\phi(x) =|x|$ for $|x|>1$ and $\phi(x)\leq |x|$  for $|x|\leq 1$,  define 
	\begin{equation}\label{Vp}
		V (x)=  1+ \phi^p(x), \quad p\in (0,1), \quad x\in \rd.
	\end{equation}
	This function $V(x)$ is the Lyapunov function which we will  use in order to show \eqref{Lyap2}.

	For $x \in \rd$ we denote by $e_x := \frac{x}{\|x\|}$ its unit vector, and by $\gamma_{x,u}$ the cosine of the angle between  the  two unit vectors:
	\begin{equation}\label{gxu}
		\gamma_{x,u}:= \cos \angle (e_x,e_u). 
	\end{equation}
	Let
	\begin{equation}\label{AB}
		\Bb(r):= \{u:\, |k(x,u)|\leq r \}, \quad \Aa(r,R):= \{u: \, r<|k(x,u)\leq R\}. 
	\end{equation}

	Fix $\lambda \in (0,1)$. For $V(x)$ as above split $\LL V$ into the ``drift'' and ``kernel'' components:
	\begin{equation}\label{split1}
		\begin{split}
			\LL V(x) & = \nabla V(x)\cdot \Big(\ell(x) + \int_{\Aa(1, |x|)}  k(x,u)\nu(du) \Big) \\
			& \quad + \int_{\rd} \left(V(x+k(x,u))- V(x) - \nabla V(x) \cdot k(x,u)\I_{|k(x,u)|\leq |x|}\right) \nu(du) \\
			& =: \LL^{drift} V(x)+ \LL^{ ker} V(x). 
		\end{split}
	\end{equation}
	It is easy to see that  for $ |x|$ big enough we have
	\begin{equation*}
		\begin{split}
			\LL^{drift} V(x) 
			&= p |x|^{p-1}\Big(|\ell(x)|\gamma_{x,\ell(x)} + \int\limits_{\Aa(1, |x|)} |k(x,u)| \gamma_{x,k(x,u)}\,\nu(du)\Big).  
		\end{split}
	\end{equation*}
	Denote
	\begin{equation*}
		\phi^{drift}(x) := |x|^{p-1}\max{\Big(|\ell(x)|,\int_{\Aa(1, |x|)} |k(x,u)\,\nu(du)\Big)}.
	\end{equation*}
	Clearly, there exists $\zeta \in [-2,2]$, such that
	\begin{equation}\label{drift}
		\LL^{drift} V(x) \leq -p\zeta \phi^{drift}(x).
	\end{equation}
	Next, split the ``kernel'' part $\LL^{ker}$ into the ``ball'' and ``tail'' components: 
	\begin{equation}\label{Lsplit}
		\begin{split}
			\LL^{ker}V(x)
			&= \int\limits_{\Bb(|x|^\lambda)} \left(V(x+k(x,u)) - V(x) - \nabla V(x) k(x,u) \right) \nu(du) \\
			&\qquad + \int\limits_{\Bb^c(|x|^\lambda)} \left(V(x+k(x,u)) - V(x)  - \nabla V(x) k(x,u) \I_{\Aa(|x|^\lambda, |x|)}\right)\nu(du)\\
			&=: \LL^{ball}V(x)+\LL^{tail}V(x).
		\end{split}
	\end{equation}
	Finally, split the ``tail'' part
	\begin{equation}\label{tailsplit}
		\begin{split}
			\LL^{tail}V(x)
			&= \int\limits_{\Aa(|x|^\lambda,|x|)} \left(V(x+k(x,u)) - V(x) - \nabla V(x)\cdot k(x,u) \right) \nu(du) \\
			&\qquad + \int\limits_{\Bb^c(|x|)} \left(V(x+k(x,u)) - V(x) \right)\nu(du)\\
			&=: \LL^{big}V(x)+\LL^{large}V(x).
		\end{split}
	\end{equation}

	To compare the impact of each component on $\LL$, consider the respective ``growth'' functions:
	\begin{equation}\label{phiPM}
		\phi^{ball}_+(x) := |x|^{p-2} \int\limits_{\Bb(|x|^\lambda)}|k(x,u)|^2 \, \nu(du) \quad \text{and} \quad \phi^{ball}_-(x) :=  |x|^{p-2} \int\limits_{\Bb(|x|^\lambda)}|k(x,u)|^2 \gamma_{ x,k(x,u)}^2\, \nu(du),
	\end{equation}
	\begin{equation}\label{phiBig}
		\phi^{big}(x) :=\frac{p}{2}|x|^{p-1} \int\limits_{\Aa(|x|^{\lambda},|x|)} |k(x,u)| \, \nu(du),
	\end{equation}
	
	\begin{equation}\label{phiLarge}
		\phi^{large}(x) := \int\limits_{\Bb^c(|x|)} |k(x,u)|^p \, \nu(du).
	\end{equation}

	\begin{lemma}\label{lem1}
		We have
		\begin{equation}\label{L30}
			\LL V(x)
			\leq  -p\zeta \phi^{drift}(x)+\frac{p}{4}(1+o(1))(\phi^{ball}_+(x)+(p-2)\phi^{ball}_-(x)) +\phi^{big}(x)+\phi^{large}(x),
		\end{equation}
		where $|x|\to \infty, \,\zeta \in [-2,2]$. 
	\end{lemma}
	
	To illustrate relations between  different components we define several constants. Let
	\begin{equation}\label{S}
		C^{ball}:= 	\limsup_{|x| \to \infty} \frac{\phi^{ball}_+(x)}{\phi^{ball}_-(x)}
	\end{equation}
	If  $C^{ball}<\infty$, then  the measure $\nu(\cdot)$ is ``thick'' in the direction of $k(x,u)$ or in the opposite one,   meaning that  $\phi^{ball}_-$ is faster than $\phi^{ball}_+$, which is crucial for ``ball''-induced ergodicity. Note that by definition $C^{ball}\geq 1$. 
	
	Next, 
	\begin{equation*}
		C^{big}:= \limsup_{|x| \to \infty} \frac{\phi^{big}(x)}{\phi^{drift}(x)}, \quad C^{large}:=  \limsup_{|x| \to \infty} \frac{\phi^{large}(x)}{\phi^{drift}(x)}. 
	\end{equation*}
	Finite $C^{big} \vee C^{large}$ assures that  the ``drift'' is faster, than the  ``tail''. Depending on the  relation between  ``drift'' and ``ball'' parts,  we can expect here either ``drift''-- or ``ball''--induced ergodicity.
	
	Similarly, we define 
	\begin{equation*}
		\tilde{C}^{drift}_+:= \limsup_{|x|\to\infty} \frac{\phi^{ball}_+(x)}{\phi^{drift}(x)}, \quad 
		\tilde{C}^{drift}_-:= \liminf_{|x|\to \infty} \frac{\phi^{ball}_-(x)}{\phi^{drift}(x)}.
	\end{equation*}
	
	Finally, depending whether $\zeta$ is positive or negative we define respectively
	\begin{equation*}
		\tilde{C}^{ball}_+:= \liminf_{|x|\to\infty} \frac{\phi^{drift}(x)}{\phi^{ball}_-(x)}, \quad 
		\tilde{C}^{ball}_-:= \limsup_{|x|\to \infty} \frac{\phi^{drift}(x)}{\phi^{ball}_-(x)}. 
	\end{equation*}
	
	In order to make our results easily applicable, we introduce some    assumptions on  lower bounds of  $\phi^{drift}$ and $\phi^{ball}$. Namely, suppose that there exist positive constants $K^{drift}, K^{ball}$ and $\gamma^{drift}, \gamma^{ball}$  such that for $ |x|\gg 1$
	\begin{equation}\label{gamdrift}
		\phi^{drift}(x)\geq  K^{drift} |x|^{\gamma^{drift}}, 
	\end{equation}
	\begin{equation}\label{gamball}
		\phi^{ball}_{-}(x)\geq  K^{ball} |x|^{\gamma^{ball}}.
	\end{equation}
	
	We introduce the conditions in such a way that one can distinguish between the ``drift-induced'' and ``noise-induced'' ergodicity, as well as  determine the part of the kernel, responsible for  it. 
	
	\begin{theorem}\label{T2}
		Case 1. (\textit{``Drift-induced'' ergodicity}). Suppose that  \eqref{drift}  holds with  $\zeta \in (0,2]$, and $C^{big}\vee C^{large}$, $C^{drift}_+\vee C^{drift}_-$  are  finite. Then \eqref{Lyap2} holds with $f(u) =   Ku^{\gamma^{drift}/p}, K > 0$, if 
		\begin{equation}\label{LI1}
			- p\zeta + \frac{p}{4}\left(\tilde{C}^{drift}_+ +(p-2)\tilde{C}^{drift}_-\right) + C^{big} + C^{large} < 0.
		\end{equation}			
		
		Case 2. (\textit{``Ball-induced'' ergodicity}). Let either $C^{big}\vee C^{large} = \infty$ or $C^{big}\vee C^{large} < \infty$,  but  
		$p<C^{big}+C^{large}$.  Denote
		\begin{equation*}
			\tilde{C}^{ball} := \limsup_{|x|\to\infty} \frac{\max{(\phi^{big}(x), \, \phi^{large}(x))}}{\phi^{ball}_-(x)}. 
		\end{equation*} 							 
		Then \eqref{Lyap2} holds with $f(u) =  Ku^{\gamma^{ball}/p}, K \!>\!0$, if $C^{ball}, \tilde{C}^{ball} < \infty$, and one of the following inequalities is satisfied:
		\begin{enumerate}
			\item[a)] 
			\begin{equation}\label{case2a}
				-p\zeta \tilde{C}^{ball}_+ + \frac{p}{4}(C^{ball} + p-2)+ 2\tilde{C}^{ball}<0, \qquad \zeta >0;
			\end{equation}
			
			\item[b)] 
			\begin{equation}\label{case2b}
				-p\zeta \tilde{C}^{ball}_- + \frac{p}{4}(C^{ball} + p-2)+ 2\tilde{C}^{ball}<0, \qquad \zeta \leq 0.
			\end{equation}
			
		\end{enumerate}			
		
	\end{theorem}

	\section{Proofs}\label{Proofs}

	\begin{proof}[Proof of Lemma~\ref{lem1}]
		We begin the proof with the behaviour of the  ``ball''  part. 
		Note that for $\theta \in  [0,1]$ and  $|k(x,u)|\leq |x|^\lambda$  we get by the triangle inequality
		\begin{equation}\label{triangle1}
			|x|\cdot |1-|x|^{\lambda -1}|\leq 		|x+\theta k(x,u)| \leq |x|(1+|x|^{\lambda-1}).
		\end{equation}
		Then,  since in $\LL^{ball}$  we have   $|k(x,u)|\leq |x|^\lambda$, $\lambda \in (0,1)$, for $|x|$ big enough we have 
		$$
		V(x +k(x,u))=1+|x+k(x,u)|^p.
		$$
		
		Therefore, by Taylor's formula, 
		\begin{equation*}
			V(x+k(x,u)) - V(x) - \nabla V(x)\cdot k(x,u) = \frac12 \int\limits_0^1 (1-\theta) k'(x,u)\cdot D^2V(x+\theta k(x,u))\cdot k(x,u) \,d\theta ,
		\end{equation*}
		where
		\begin{equation*}\label{DV}
			D^2V(x) = \left\{\partial_{ij}^2 V(x)\right\}_{i,j=1}^d=p|x|^{p-2}
			\begin{pmatrix}
				1+(p-2)\frac{x_1^2}{|x|^2} & (p-2)\frac{x_1 x_2}{|x|^2} & \cdots& (p-2)\frac{x_1 x_d}{|x|^2} \\ 
				(p-2)\frac{x_1 x_2}{|x|^2} & 1+(p-2)\frac{x_2^2}{|x|^2} & \cdots & (p-2)\frac{x_2 x_d}{|x|^2}\\
				\cdots & \cdots & \ddots & \cdots\\
				(p-2)\frac{x_1 x_d}{|x|^2} & (p-2)\frac{x_2 x_d}{|x|^2} & \cdots &  1+(p-2)\frac{x_d^2}{|x|^2}
			\end{pmatrix}.
		\end{equation*}
		Then 
		\begin{align*}
			\LL^{ball} V(x) 
			&= \frac{1}{2}  \int\limits_{\Bb(|x|^\lambda)}\left[\int\limits_0^1 (1-\theta) k'(x,u) \cdot D^2V\left(x+\theta k(x,u)\right) \cdot k(x,u) \, d\theta \right] \, \nu(du)\\
			&=\frac{p}{2} \int\limits_{\Bb(|x|^\lambda)} \int\limits_0^1 (1-\theta)|x+\theta k(x,u)|^{p-2}\left(\sum_{i=1}^{d}{\left(1+(p-2)\frac{(x_i+\theta k_i(x,u))^2}{|x+\theta k(x,u)|^2}\right)k_i^2(x,u)}\right.\\
			&\qquad \qquad \left.+(p-2)\sum_{i,j=1, i\neq j}^{d}{k_i(x,u) k_j(x,u)\frac{\left(x_i+ \theta k_i(x,u)\right) \left(x_j+\theta k_j(x,u)\right)}{|x+\theta k(x,u)|^2}}\right) d\theta \nu(du).
		\end{align*}
		Now   using the sum of squares  formula for $k_i(x,u)(x_i + \theta k_i(x,u))$, we derive 
		\begin{align*}
			\LL^{ball} V(\!x\!)
			&= \frac{p}{2} \!\int\limits_{\Bb(|x|^\lambda)} \!\int\limits_0^1 \frac{1-\theta}{|x+\theta k(x,\!u)|^{2-p}} {\left(|k(x,u)|^2 + (p-2)\! \left(\!k(x,u)\frac{x+\theta k(x,u)}{|x+\theta k(x,u)|}\!\right)^2\right)}\!d\theta\nu(du)\\
			&= \frac{p}{2} \int\limits_{\Bb(|x|^\lambda)}\left[\int\limits_0^1 \frac{(1-\theta)|k(x,u)|^2}{|x+\theta k(x,u)|^{2-p}}\left(1 + (p-2)\left(e_{k(x,u)}\cdot e_{x+\theta k(x,u)}\right)^2\right) d\theta\right] \nu(du).
		\end{align*}
		
		Using  \eqref{triangle1} we get 
		\begin{equation}\label{ball1}
			\int\limits_0^1 \frac{1-\theta}{|x+\theta k(x,u)|^{2-p}} \,d\theta \leq  \frac{|x|^{p-2}}{2\left(1-|x|^{\lambda-1}\right)^{2-p}}
		\end{equation}	
		and 
		\begin{equation}\label{ball2}
			(p-2)\int\limits_0^1 \frac{1-\theta}{|x+\theta k(x,u)|^{2-p}} \, d\theta  
			\leq \frac{p-2}{2}\frac{|x|^{p-2}}{\left(1+ |x|^{\lambda-1}\right)^{2-p}}.
		\end{equation}
		
		Finally,  transforming  the scalar  product of the vectors  $e_{k(x,u)}$, $e_{x+\theta k(x,u)}$  and applying the estimate \eqref{triangle1}, we get
		\begin{equation}\label{dotp}
			\begin{split}
				|e_{k(x,u)}\cdot e_{x+\theta k(x,u)} |^2
				&=  \left|  e_{k(x,u)}\cdot e_{x}\frac{|x|}{|x+\theta k(x,u)|} +   \frac{\theta |k(x,u)| }{|x+\theta k(x,u)|} \right|^2
				\\& 
				\geq \gamma_{x,k(x,u)}^2 \left(\frac{|x|}{|x|+ |k(x,u)|}\right)^2 - \frac{2\theta |x|  \cdot |k(x,u)| }{(|x|-|k(x,u)|)^2}
				\\& \geq 
				\frac{\gamma_{x,k(x,u)}^2 }{(1+|x|^{\lambda-1})^2}- 
				\frac{2 |x|^\lambda}{|x|(1-|x|^{\lambda-1})^2}.
			\end{split}
		\end{equation}
		
		Combining  \eqref{ball1}, \eqref{ball2} and \eqref{dotp}, we derive
		\begin{align*}
			\LL^{ball}(x)
			&\leq \frac{p}{4} |x|^{p-2} \int\limits_{\Bb(|x|^\lambda)}|k(x,u)|^2 \left(C_1(x)+ (p-2)C_2(x) \gamma^2_{x,  k(x,u)}\right)\nu(du) \\
			&= \frac{p}{4} C_1(x) \phi^{ball}_+(x) +  \frac{p(p-2)}{4} C_2(x) \phi^{ball}_-(x),
		\end{align*}
		where  
		\begin{equation}\label{C12}
			C_1(x) =\frac{1}{\left(1-|x|^{\lambda-1}\right)^{2-p}}  +  \frac{2(2-p) |x|^{\lambda-1}}{(1+|x|^{\lambda -1})^{2-p}(1-|x|^{\lambda -1})^2}, \quad 
			C_2(x) = \frac{1}{\left(1+ |x|^{\lambda-1}\right)^{4-p}}.
		\end{equation} 
		Note that both $C_1(x)$ and $C_2(x)$ tend to 1 as $|x| \to \infty$, hence 
		\begin{equation}\label{Lball}
			\LL^{ball} V(x)
			\leq  \frac{p}{4}(1+o(1))(\phi^{ball}_+(x)+(p-2)\phi^{ball}_-(x)). 
		\end{equation}
		
		Consider now the ``tail''. For the big jumps
		\begin{align*}
			\LL^{big} V(x)
			&\leq \int\limits_{\Aa(|x|^\lambda, |x|)} (|x+k(x,u)|^p-|x|^p-p|x|^{p-1}\cdot |k(x,u)|\gamma_{ x,k(x,u)})\, \nu(du) \\
			&= |x|^p \int\limits_{\Aa(|x|^\lambda, |x|)} \Big(\Big(1 + 2\gamma_{ x,k(x,u)} \frac{|k(x,u)|}{|x|} + \Big(\frac{|k(x,u)|}{|x|} \Big)^2\Big)^{\frac{p}{2}}  -1\Big)\, \nu(du) \\
			&-p|x|^{p-1} \int\limits_{\Aa(|x|^\lambda, |x|)}|k(x,u)|\gamma_{ x,k(x,u)}\, \nu(du) =: J_1(x)+J_2(x).
		\end{align*}
		Since on $\Aa(|x|^\lambda, |x|)$ we have $|k(x,u)| \leq |x|$, then
		\begin{equation*}
			\begin{split}
				J_1(x)
				&\leq |x|^p \int\limits_{\Aa(|x|^\lambda, |x|)} \Big(\Big(1 + (1+2\gamma_{ x,k(x,u)})\frac{|k(x,u)|}{|x|} \Big)^{\frac{p}{2}}  -1\Big) \nu(du) \\
				&\leq \frac{p}{2}|x|^{p-1} \int\limits_{\Aa(|x|^\lambda, |x|)} (1+2\gamma_{ x,k(x,u)})|k(x,u)|\,\nu(du).
			\end{split}
		\end{equation*}	
		Combined with $J_2$, the above estimate yields 
		\begin{equation}\label{big}
			\LL^{big} V(x) \leq \frac{p}{2}|x|^{p-1} \int\limits_{\Aa(|x|^\lambda, |x|)} |k(x,u)| \nu(du) = \phi^{big}(x). 
		\end{equation}
		
		For the large jumps we simply have
		\begin{equation}\label{large}
			\LL^{large}V(x) 
			\leq  \int\limits_{\Bb^c(|x|)}|k(x,u)|^p \nu(du)
			= \phi^{large}(x).
		\end{equation}

		Combining \eqref{drift}, \eqref{Lball}, \eqref{big} and \eqref{large}, we derive  the statement of the lemma: 
		\begin{equation*}
			\LL V(x)
			\leq  -p \zeta \phi^{drift}(x)+\frac{p}{4}(1+o(1))(\phi^{ball}_+(x)+(p-2)\phi^{ball}_-(x)) +\phi^{big}(x)+\phi^{large}(x), \, |x| \to \infty. \qedhere 
		\end{equation*}
	\end{proof}

	\begin{proof}[Proof of Theorem~\ref{T2}] 
		\textit{Case 1}. Consider first the ``drift-induced'' ergodicity. If  $\zeta \in  (0, 2]$, the drift is negative. Divide both sides of \eqref{L30} by $\phi^{drift}(x)$ and take $\limsup$ as $|x|\to \infty$:
		\begin{equation}\label{Ld}
			\limsup_{|x| \to \infty} \frac{\LL V(x)}{\phi^{drift}(x)}	 \leq - p\zeta + \frac{p}{4}\left(\tilde{C}^{drift}_+ +(p-2)\tilde{C}^{drift}_-\right)+ C^{big} + C^{large}.
		\end{equation}
		Since both $C^{big}$ and $C^{large}$ are positive, we need them to be finite to ensure that the drift is faster than the tail. Then the Lyapunov-type condition \eqref{Lyap2} is satisfied, if the right-hand side of \eqref{Ld} is negative.

		\textit{Case 2}.  The ``ball induced''  ergodicity may appear,  if either the tail is faster than the drift ($C^{big}\vee C^{large} = \infty$) or the drift (if negative) is fast enough, but still it cannot compensate the positivity of the tail, i.e.  $C^{big}\vee C^{large} < \infty$, but  $p <C^{big}+C^{large}$. In this scenario in order to have  \eqref{Lyap2},  the ball part  has to be negative and fast enough, i.e $C^{ball}<2-p$ and $\tilde{C}^{ball}<\infty$. Dividing both sides of \eqref{L30} by $\phi^{ball}_-(x)$, taking $\limsup$ as $|x|\to \infty$, we see that the sufficient condition for  \eqref{Lyap2}  is either \eqref{case2a}, if $\zeta>0$, or \eqref{case2b}, if $\zeta<0$. 
		
		In both cases the choice of $f$ in \eqref{Lyap2} follows from \eqref{gamdrift} and \eqref{gamball}. \qedhere 
	\end{proof}
	
	\section{Examples}\label{Examples}
	
	\begin{example}
		$d=1$. 	Consider a one-sided $\alpha$-stable  L\'evy measure  $\nu(du) = u^{-1-\alpha}\I_{u>0}du$ (we assume  for simplicity that the normalizing constant equals  1). Take $k(x,u) =ub(x)$, where $b(\cdot)$ is Lipschitz and  $ b(x)= |x|^{\beta }\sign(-x)(1+o(1))$ as $x\to\pm\infty$  with  $\beta \in (0,1)$;  $\ell(x) = |x|^\eta$, $ \eta \in (0,\alpha \beta)$. Condition \eqref{A1}  is clearly satisfied. In order to see that  \eqref{A2} holds true as well, one can use the fact  that the process $X$ possesses a strictly positive transition probability density, see, for example, \cite{FG06}.
		
		To make calculations a bit more transparent let us introduce an auxiliary function
		\begin{equation}\label{phiAnn}
			\phi^{ann}_\tau(x)
			:= |x|^{p-1} \int\limits_{\Aa(|x|^\tau,|x|)} |k(x,u)| \nu(du), \quad \tau \in \{0,\lambda\}. 
		\end{equation} 
		
		Depending on the value of $\alpha$ we have
		
		\begin{equation*}
			\phi^{ann}_\tau(x)  \asymp 
			\begin{cases} 
				(\alpha-1)^{-1}|x|^{ p-1+\alpha\beta+\tau (1-\alpha)}, & \alpha>1,\\
				(1-\alpha)^{-1}|x|^{ p-\alpha+\alpha\beta}, & \alpha<1, \\
				(1-\tau)|x|^{ p-1 + \beta }\ln |x|,	& \alpha =1. 
			\end{cases}
		\end{equation*}
		Then
		\begin{equation*}
			\phi^{drift}(x)  = \phi^{ann}_0(x) \quad \text{and} \quad\phi^{big}(x)  = \frac{p}{2}\phi^{ann}_\lambda(x).
		\end{equation*}
		
		Next
		$$
		\phi^{ball}_+(x) = \phi^{ball}_-(x) =\frac{1}{2-\alpha} |x|^{p-2 + \alpha\beta +\lambda(2-\alpha)} (1+ o(1)), 
		$$
		and	
		$$
		\phi^{large}(x) = \frac{1}{\alpha - p} |x|^{p-\alpha+ \alpha \beta} (1+o(1)). 
		$$

		Note that $\tilde{C}^{drift}_\pm = 0$ for all $\alpha \in (0,2)$, which means that the  ``ball'' is slower than the  ``drift'', and that we are in the settings of the case of the ``drift-induced'' ergodicity,  but the ``kernel part'' of the drift dominates due to the assumption $\eta<\alpha \beta$.  Moreover $\zeta = 1$ for all $\alpha \in (0,2)$.

		If $\alpha >1$, $C^{big} = C^{large} = 0$.  Then \eqref{Lyap2} holds with $f(u) = Cu^{\frac{p-1+\alpha \beta}{p}}$, $C \leq p$.  
		
		If $\alpha <1$, $C^{big} = \frac{p}{2}$, $C^{large}=\frac{1-\alpha}{\alpha-p}$. Then \eqref{Lyap2} holds with $f(u) = Cu^{\frac{p-\alpha+\alpha \beta}{p}}$, $C \leq \frac{p}{2}-\frac{1-\alpha}{\alpha-p}$, provided $p^2-p\alpha+2(1-\alpha)<0$. 
		
		If $\alpha =1$, $C^{big} = \frac{p}{2}(1-\lambda)$, $C^{large}=0$.  Then \eqref{Lyap2} holds with $f(u) = Cu^{\frac{p-1+\beta}{p}}\ln{u}$, $C \leq \frac{p}{2}(1+\lambda)$. 
		
		It is a known result (see, for example,  \cite{BS15}), that in this case the limit probability measure  is trivial.
	\end{example}

	\begin{example}
		Let $d=1$. Consider a kernel, such that $\nu\{\pm K \}\!=\!1$ for some $K>0$. Take again $k(x,u) = |x|^\beta u$, $\beta\in (0, \lambda)$,  and $\ell(x)$  such that $|\ell(x)| = |x|^\eta, \eta > 1-p$ and  $x\ell(x) <0$. In this case, 
		\begin{equation*}
			\phi^{ball}_+ (x)= \phi^{ball}_- (x)  = 2K^2|x|^{2\beta}|x|^{p-2}
		\end{equation*}
		and
		\begin{equation*}
			\phi^{drift}(x) =  - |x|^{p-2} |x\ell(x)|= - |x|^{p-1+\eta}.
		\end{equation*}
		Moreover, $C^{ball}=1$, and for $|x|$ big enough  $\phi^{big}(x)$, $\phi^{large}(x)$ tend to zero. Therefore,
		\begin{equation}\label{L20}
			\LL V(x) 
			\leq -p|x|^{p-1}
			\left(|x|^\eta + \frac{(1-p)K^2}{2} |x|^{2\beta-1} \right).
		\end{equation}
		Then for $\eta \leq 2\beta -1$ the Lyapunov inequality \eqref{Lyap2} holds with  $f(u) = C u^{\frac{p-2+2\beta}{p}}$,  $C \leq \frac{K^2}{2}(1-p)$. Otherwise \eqref{Lyap2} holds with $f(u) = C u^{\frac{p-1+\eta}{p}}$, $C \leq p$.
		
		The limit probability measure is again trivial.
	\end{example}
	
	\begin{example}
		Consider now a multidimensional case, but with zero drift part $\ell$. Take a multiplicative $k(x,u) = \Phi(x)u$, where $\Phi \colon \rd\mapsto \rd\times \rd$ is a matrix, whose elements satisfy  the Lipschitz and linear growth conditions from \cite{GS82}. We assume that $\Phi(x) =|x|^\beta U(x)$, where $U(x)$ is a unitary matrix  and $\beta > 1-p/2$.  Assume also that $\nu(\cdot)$ is symmetric and concentrated entirely on the sphere in $\rd$ centered at $0$  and of radius $R$. Then the tail part $\LL^{tail}$ is simply absent, and the annulus part turns into zero due to the symmetry of $\nu(\cdot)$. 
		
		Since matrix $U$ is unitary, we have $|\Phi(x) u | = |x|^{\beta}|u|$, and thus
		\begin{align*}\label{phiPhi}
			\phi^{ball}_-(x)&
			= |x|^{p-2} |x|^{2\beta}\int\limits_{ |u| =R} |u|^2 \gamma_{\Phi(x)u, x}^2 \, \nu(du) \\
			&= |x|^{p-2} |x|^{2\beta} R^2 \int\limits_{|u| =R} \gamma_{\Phi(x)u, x}^2 \, \nu(du) \notag.
		\end{align*}
		Therefore,
		\begin{equation*}\label{SPhi}
			C^{ball}
			= \frac{\nu(u: |u|=R)}{\int\limits_{|u|=R} \gamma_{\Phi(x)u, x}^2\, \nu(du)}.
		\end{equation*}
		If we chose $p<\min(2-C^{ball},1)$, we can achieve that \eqref{Lyap2} holds true   with 
		$$
		f(u) = C u^{\frac{p-2+2\beta}{p}}, \quad  \text{where}\quad C  \leq \frac{p}{4} R^2 (2-p -C^{ball}). 
		$$
	\end{example}

	\subsection*{Acknowledgement.} 
	V. Knopova was supported by  the internal Grant 22BF038-01 of the Research Department of Taras Shevchenko National University of Kyiv; V. Knopova also thanks the V.M.Glushkov Institute of Cybernetics NAS of Ukraine for the hospitality during her probation time in May-June 2024. 
	
	This research was partly funded by the ``DAAD Ostpartnerschaften Programme'' of the TU Dresden. We are grateful to Faculty of Mathematics at TU Dresden and Prof.\ Ren\'e Schilling for providing excellent working conditions.
	
	The authors also wish to express their sincere gratitude to Prof.\ Ren\'e Schilling for his insightful advice.
	
	The authors thank Prof. O.Kulyk for his comments, which helped to improve the paper. The authors also thank the anonymous reviewer for their valuable advice and comments.
	
	\nocite{*}
	
	\bibliographystyle{unsrtnat}
	
	\bibliography{UnbddArxiv2025}

\end{document}